\documentclass{amsart}
\usepackage[cp1250]{inputenc}
\usepackage[]{ tikz, hyperref}

\usepackage{geometry}
\usepackage{graphicx} % Required for inserting images

\usetikzlibrary{automata, positioning}

\newtheorem{theorem}{Theorem}[section]
\newtheorem{proposition}[theorem]{Proposition}
\newtheorem{lemma}[theorem]{Lemma}

\newtheorem{definition}[theorem]{Definition}
\newtheorem{example}[theorem]{Example}
\newtheorem{remark}[theorem]{Remark}

\numberwithin{equation}{section}

\title{Three combinatorial identities of the Rogers-Ramanujan type}
\author{Stefano Capparelli, Arne Meurman, and Mirko Primc}
\begin{document}

\maketitle

\section*{Abstract}
The purpose of the present paper is to establish three
combinatorial Rogers--Ramanujan type identities. Although their
motivation comes from the representation theory of affine Lie algebra of type $A_2^{(2)}$, our
proofs are entirely combinatorial and analytic. They follow the
automata-theoretic approach introduced by Takigiku and Tsuchioka in
their proof of Nandi's conjecture and adapt it to the present setting.
%%%%%%%%%%%%%%%

\section{Introduction}

\noindent
The Rogers--Ramanujan identities have inspired a remarkable body of
mathematics. Their partition-theoretic interpretation established deep
connections between infinite products, integer partitions, and
$q$-series, while later work of Lepowsky--Milne and Lepowsky--Wilson
revealed equally striking connections with the representation theory of
affine Lie algebras. In this setting, specialized characters of standard
modules frequently admit two complementary descriptions. On the one
hand, the Weyl--Kac character formula, together with Lepowsky's
numerator formula, yields explicit infinite-product expressions. On the
other hand, one seeks combinatorial realizations of these products as
generating functions for partitions satisfying suitable local difference
conditions. Many celebrated partition identities have emerged from this
interaction between representation theory and combinatorics.

The present work was originally motivated by the study of the three
level-five standard modules of the twisted affine Lie algebra
$A_2^{(2)}$. The principal specializations of their characters suggest
three infinite products of Rogers--Ramanujan type, naturally leading to
the search for corresponding partition identities. The first of these
products
$$
\prod_{\underset{n\equiv \pm 2,\pm 3,\pm 4, \pm 5, \pmod{16}}{n\in\mathbb{N}}}(1-q^n)^{-1},
$$
associated with the highest weight $5\Lambda_0$, was
investigated through extensive computer experimentation, eventually
leading to a family of local difference conditions (see 2.2 below) whose generating
function agreed with the product to all computed orders. Although these
conditions were initially sought as a possible combinatorial
realization of the corresponding module, subsequent evidence suggested
that such an interpretation was unlikely to hold. Nevertheless, the
resulting partition identity proved to be of independent combinatorial
interest. The same investigation ultimately led to analogous identities
corresponding to the remaining two level-five products.

The purpose of the present paper is to establish these three
combinatorial Rogers--Ramanujan type identities, i.e., Theorems \ref{theorem1}, \ref{theorem2} and \ref{theorem3}. Although their
motivation comes from the representation theory of $A_2^{(2)}$, our
proofs are entirely combinatorial and analytic. They follow the
automata-theoretic approach introduced by Takigiku and Tsuchioka in
their proof of Nandi's conjecture \cite{TT19} and adapt it to the present setting.
This extends the methods of \cite{A76} Chapter 8.

Our approach begins by reformulating the partition conditions in terms
of frequency vectors. This converts the problem into the study of finite
local constraints on words over the alphabet
\[
\Sigma=\{0,1,2\}.
\]
The local nature of these constraints leads naturally to deterministic
finite automata recognizing admissible frequency vectors. The
continuation generating functions associated with the states of these
automata satisfy finite systems of $q$-difference equations. After
suitable reductions, these systems are solved explicitly by identifying
appropriate double-sum solutions and proving that they are uniquely
determined by the corresponding functional equations.
To obtain the product sides we use three closely related identities from \cite{AEH20}, 
where they arise from characters of three modules for the Virasoro Lie algebra.
Since the automata-theoretic method is likely to be unfamiliar to many
readers, we develop it in full detail for the first identity. 

The paper is organized as follows. Sections~2--4 develop the
automata-theoretic framework for the first identity, beginning with the
translation of the partition conditions into local constraints on
frequency vectors and culminating in the construction of the associated
deterministic finite automaton. The subsequent sections derive and solve
the corresponding system of $q$-difference equations and establish the
first Rogers--Ramanujan type identity. The final sections describe the
necessary modifications for the remaining two identities and apply the
same method to complete their proofs.
The three partition families exhibit an evident parallelism. Although
their initial conditions differ, reflecting the distinct level-five
modules that motivated them, their local difference conditions have the
same general structure. Consequently, the associated automata, the
functional equations satisfied by their continuation generating
functions, and the subsequent proofs are closely analogous. We therefore
develop the automata-theoretic construction in full detail only for the
first identity and treat the remaining two identities more
economically, emphasizing only the features that differ from the
prototype. 

After submission to arXiv of the first version of this paper, we have learned that three similar combinatorial identities have been proved in \cite{T22}.

\section{Frequency vectors and local conditions}\label{frequency_vectors}

\noindent
The partition identities studied in this paper are most naturally formulated
in terms of frequency vectors rather than partitions themselves.
This representation makes the local nature of the difference conditions
transparent and provides the language in which the subsequent
automata-theoretic construction will be developed.

\subsection{Partitions and frequency vectors}

A \emph{partition} of a nonnegative integer $n$ is a finite nonincreasing
sequence of positive integers
\[
\lambda=(\lambda_1,\lambda_2,\ldots,\lambda_r),
\]
whose sum is $n$.
The integer $\lambda_i$ is called a \emph{part} of the partition.

For each positive integer $j$, let
\[
f_j=f_j(\lambda)
\]
denote the multiplicity of the part $j$ in $\lambda$.
The infinite sequence

\[
\mathbf f(\lambda)
=(f_1,f_2,f_3,\ldots)
\]
is called the \emph{frequency vector} of $\lambda$.

Conversely, every sequence of nonnegative integers having finite support
determines a unique partition.
Thus partitions and finite-support frequency vectors are in bijective
correspondence.

Having reformulated the conditions in terms of frequency vectors, we adopt this description for the remainder of the paper.

\begin{example}
The partition
\[
(8,8,7,5,5,3,2,2)
\]
has frequency vector
\[
(f_1,f_2,f_3,f_4,f_5,f_6,f_7,f_8,f_9, \dots)
=
(0,2,1,0,2,0,1,2,0, \dots).
\]
\end{example}

\begin{definition}
  A partition $\mu$ is called a \emph{subpartition} of a
  partition $\lambda$ if, for every positive integer $j$, the multiplicity of $j$ in $\mu$ does not exceed its multiplicity in   $\lambda$, (cf. \cite{N14}, \cite{AEH20}).
\end{definition}

\subsection{Difference conditions}
The first family of partitions considered in this paper is denoted by $\mathcal{C}$ and defined by the following initial condition

(IC1): No part is equal to $1$.

\noindent
and by ten local difference conditions (DC), originally discovered during the computational search described in the Introduction.

\begin{description}

\item[(D1)]
No part occurs three or more times.

\item[(D2)]
No partition has a subpartition $(k,k,k-1)$.

\item[(D3)]
No partition has a subpartition $(k,k-1,k-1)$.

\item[(D4)]
No partition has a subpartition $(k,k-1,k-2)$.

\item[(D5)]
No partition has a subpartition $(k,k,k-2)$.

\item[(D6)]
No partition has a subpartition $(k,k-1,k-3)$.

\item[(D7)]
No partition has a subpartition $(k,k,k-3,k-3)$.

\item[(D8)]
No partition has a subpartition $(k,k-1,k-4,k-4)$.

\item[(D9)]
No partition has a subpartition $(k,k,k-3,k-4)$.

\item[(D10)]
No partition has a subpartition $(k,k,k-3,k-5,k-5)$.

\end{description}
We also say that a partition in $\mathcal{C}$ is an {\em admissible partition}.

\subsection{Immediate consequences}

Several simple consequences of conditions (IC1) and (D1)--(D10) will be used
throughout the paper.

\begin{remark}
Condition (D1) immediately implies that
\[
f_j\in\{0,1,2\}
\]
for every positive integer \(j\). Thus every frequency vector may be
viewed as an infinite word over the alphabet
\[
\Sigma=\{0,1,2\}.
\]
\end{remark}

\begin{remark}
Each of the conditions (D1)--(D10) involves only finitely many
consecutive frequencies.
Consequently, admissibility of a partition is completely determined by local
properties of its frequency vector.
\end{remark}

 The locality of these conditions permits an equivalent reformulation in terms of forbidden finite words over the alphabet $\Sigma=\{0,1,2\}$. This reformulation is the subject of the next section.

\section{Reduction to forbidden local configurations}
\label{sec:frequency}

\subsection{Translation of the difference conditions}

In this section we reformulate the original partition conditions as local
conditions on the frequency vector.

Throughout this section, local configurations of a frequency vector are
written in increasing order of the indices. Thus, for example, the
adjacent pattern of length $s+1\geq2$
\[
(a_0, \dots, a_s)
\]
means that
\[
(f_i, \dots, f_{i+s})=(a_0, \dots, a_s)
\]
for some \(i\geq 1\).

%%%%%

The difference conditions (D1)--(D10) may be
reformulated entirely in terms of local restrictions on the associated
frequency vector. Moreover, each restriction reduces to a finite collection
of irreducible forbidden local configurations.

\begin{proposition}\label{prop:forbidden}
A partition satisfies conditions \textup{(IC1)} and \textup{(D1)}--\textup{(D10)} if and
only if its frequency vector satisfies the following conditions.

\begin{enumerate}
\item[(F1)] \(f_1=0\).
\item[(F2)] \(f_i\in\{0,1,2\}\) for every \(i\ge1\) and
\end{enumerate}
the following adjacent configurations are forbidden
\begin{enumerate}
\item[(F3)] $(1,2),\quad(2,1),\quad(2,2),$

\item[(F4)] $(1,1,1),$

\item[(F5)] $(1,0,2),\quad (2,0,2),$

\item[(F6)] $(1,0,1,1), \quad (2,0,1,1),$

\item[(F7)] $(2,0,0,2),$

\item[(F8)] $(2,0,0,1,1),$

\item[(F9)] $(1,1,0,0,2),$

\item[(F10)] $(2,0,1,0,0,2).$
\end{enumerate}
\end{proposition}

%%%%%
\begin{proof}
Condition \textup{(IC1)} is equivalent to
\[
f_1=0,
\]
since \(f_1\) is the multiplicity of the part \(1\).
Unlike the remaining conditions, this is an initial restriction on the
frequency vector rather than a translation-invariant local condition.

Condition \textup{(D1)} is equivalent to
\[
f_i\leq 2
\qquad\text{for every } i\geq 1.
\]
Since the frequencies are nonnegative integers, this is the same as
\[
f_i\in\{0,1,2\}
\qquad\text{for every } i\geq 1.
\]

Condition \textup{(D2)} forbids a subpartition of the form
\[
(k,k,k-1).
\]
In terms of the frequency vector, this means that
\[
f_k=2
\qquad\text{and}\qquad
f_{k-1}\geq 1
\]
cannot occur simultaneously.

Since, by condition \textup{(D1)}, each frequency belongs to
\(\{0,1,2\}\), this is equivalent to forbidding the adjacent frequency
patterns
\[
(f_{k-1},f_k)=(1,2)
\qquad\text{and}\qquad
(f_{k-1},f_k)=(2,2).
\]

Condition \textup{(D3)} forbids a subpartition of the form
\[
(k,k-1,k-1).
\]
In terms of the frequency vector, this means that
\[
f_{k-1}=2
\qquad\text{and}\qquad
f_k\geq1
\]
cannot occur simultaneously.

Since, by condition \textup{(D1)}, each frequency belongs to
\(\{0,1,2\}\), this is equivalent to forbidding the adjacent frequency
patterns
\[
(f_{k-1},f_k)=(2,1)
\qquad\text{and}\qquad
(f_{k-1},f_k)=(2,2).
\]

Conditions \textup{(D2)} and \textup{(D3)} together imply that no entry
equal to \(2\) may be adjacent to a positive entry. Equivalently, the
adjacent frequency patterns
\[
(1,2),\ (2,1),\ (2,2)
\]
are all forbidden.

Condition \textup{(D4)} forbids a subpartition of the form
\[
(k,k-1,k-2).
\]
In terms of the frequency vector, this means that
\[
f_{k-2}\geq 1,\qquad
f_{k-1}\geq 1,\qquad
f_k\geq 1
\]
cannot occur simultaneously.

Since the frequencies belong to \(\{0,1,2\}\), this initially forbids
all configurations
\[
(f_{k-2},f_{k-1},f_k)=(a,b,c),
\qquad a,b,c\in\{1,2\}.
\]
Conditions \textup{(D2)} and \textup{(D3)}, however, already forbid a
frequency equal to \(2\) from being adjacent to a positive frequency.
Consequently, every such configuration containing a \(2\) is already
excluded. Thus the only new forbidden local configuration arising from
\textup{(D4)} is
\[
(1,1,1).
\]

Conditions \textup{(D5)}--\textup{(D10)} are treated in a similar way.
\end{proof}
%%%%%%

The reformulation in Proposition~\ref{prop:forbidden} is the key
combinatorial step of the proof. It replaces the original partition
conditions by a finite family of local constraints on words $\Sigma^*$ over the
alphabet
\[
\Sigma=\{0,1,2\}.
\]
Consequently, admissibility of a frequency vector depends only on
finite local information. This observation makes it possible to
recognize admissible frequency vectors by means of a deterministic
finite automaton, which is constructed in the next section.

\section{Construction of the finite automaton}\label{dfaconstruction}

\noindent
In the preceding section the original partition conditions were reduced
to a finite collection of forbidden local configurations of the
associated frequency vector. Since these configurations involve only
finitely many consecutive frequencies, admissibility depends only on
finite local information. This observation naturally leads to a
deterministic finite automaton recognizing admissible frequency vectors.

Rather than introducing the automaton directly, we first determine what
information must be retained after reading a finite prefix of a
frequency vector in order to decide which continuations remain
admissible.

In Section 2 we introduced the frequency vector ${\bf f}(\lambda)$ as an infinite sequence, which always ends in an infinite sequence of zeros. Below, we shall identify this with any of its finite truncations that contain all nonzero entries. For example,
$(1,0,1,0,2,0,0)$ is identified with ${\bf f}((5,5,3,1))$.

\subsection{Finite memory}

\begin{definition}\label{defDFA}\cite{K97}
    A {\em deterministic finite automaton} (DFA) is a fivetuple $M = (Q,\Sigma,\delta,s,A)$, where $Q$ is a finite set, called the set of {\em states} of the automaton, $\Sigma$ is a finite set called the {\em alphabet}, $\delta:Q\times \Sigma\to Q$ is a map giving the {\em transitions} of $M$, $s\in Q$ is a distinguished state called the {\em start state} of $M$, and $A$ is a subset of $Q$, called the set of {\em accepting states} of $M$. 
\end{definition}

Here we have changed the notation in [TT19] for the set of accepted states from $F$ to $A$.
Let $M = (Q,\Sigma,\delta,s,A)$ be a DFA and let $w = i_1...i_k$ be a word built from elements $i_j\in\Sigma$. Then $w$ is said to be {\em accepted} by $M$ if $\delta(\dots\delta(s,i_1),\dots,i_k)\in A$. The set of words accepted by $M$ is called the {\em language}, $L(M)$, {\em accepted} by $M$.

In this section, we shall construct a DFA $M_1$ such that the language accepted by $M_1$ is the set of admissible frequency vectors. 
Let

\[
u=(f_1,f_2,\ldots,f_n)
\]
be an admissible finite frequency vector.

Suppose that a new frequency
\[
a\in\{0,1,2\}
\]
is appended to obtain

\[
(f_1,f_2,\ldots,f_n,a).
\]

The following proposition shows that only the most recent five
frequencies are relevant when deciding whether the new word remains
admissible.

\begin{proposition}\label{prop:memory}
Whether a symbol
\[
a\in\{0,1,2\}
\]
may be appended to an admissible frequency vector depends only on the
last five frequencies already read.
\end{proposition}

\begin{proof}
By Proposition~\ref{prop:forbidden}, every forbidden local configuration
has length at most six, the longest being

\[
(2,0,1,0,0,2).
\]

Appending a single frequency cannot affect the admissibility of any
configuration that does not contain the newly appended symbol.
Consequently, a new violation can occur only if the forbidden
configuration ends at the appended frequency.

Since every forbidden configuration has length at most six, only the
preceding five frequencies can participate in such a configuration.
Hence the admissibility of the appended frequency depends only on the
last five frequencies of the previously read word.
\end{proof}

Proposition~\ref{prop:memory} suggests representing the current state of
the computation by the suffix consisting of the last five frequencies.
However, distinct suffixes need not lead to different future behaviour.
Indeed, many admissible five-entry suffixes admit exactly the same
admissible continuations.

The next subsection introduces an equivalence relation on admissible
five-entry suffixes that identifies suffixes having identical future
behaviour. The resulting equivalence classes will form the states of the
automaton.

\subsection{Continuation classes}

Proposition~\ref{prop:memory} shows that the admissibility of future
extensions depends only on the current five-entry suffix. A first
approach would therefore be to construct a finite automaton having one
state for every admissible five-entry suffix.

This, however, would introduce unnecessary redundancy. Distinct
five-entry suffixes may admit exactly the same admissible future
extensions and therefore carry precisely the same information. From the
point of view of all subsequent continuations, such suffixes are
indistinguishable and may be identified.

We formalize this idea as follows. We shall use the notation $\Sigma^*$ for 
the set of words on $\Sigma$.

\begin{definition}
Let \(u\) and \(v\) be admissible five-entry suffixes.

We say that \(u\) and \(v\) belong to the same \emph{continuation class}
if, for every finite word
\[
w\in\{0,1,2\}^*,
\]
the word $uw$ obtained by appending \(w\) to \(u\) is admissible if and only
if the word $vw$ obtained by appending \(w\) to \(v\) is admissible.
\end{definition}

\begin{remark}
A continuation class is not a single five-entry suffix. Rather, it
consists of all admissible suffixes having exactly the same admissible
future extensions.

The finite automaton constructed below will have one state
corresponding to each continuation class.
\end{remark}

The continuation classes provide the natural candidates for the states
of a finite automaton. To justify this construction, two facts must be
established. First, the reachable admissible five-entry suffixes must be
partitioned into finitely many continuation classes. Secondly, it must
be shown that appending a symbol to a suffix induces a well-defined
transition between continuation classes.

We first classify the reachable suffixes.
\begin{proposition}\label{prop:classes}
The relation of belonging to the same continuation class is an
equivalence relation on the set of reachable admissible five-entry
suffixes.
Moreover, this equivalence relation has exactly eight equivalence
classes. These classes are listed in
Table~\ref{tab:classes}.
\end{proposition}
\begin{proof}
Reflexivity and symmetry are immediate from the definition.
For transitivity, suppose that \(u,v\) belong to the same continuation
class and that \(v,z\) belong to the same continuation class.
Then every finite continuation is admissible for \(u\) if and only if it
is admissible for \(v\), and likewise for \(v\) and \(z\).
Hence the same equivalence holds for \(u\) and \(z\).

Starting from the initial suffix
\[
(0,0,0,0,0),
\]
one enumerates all reachable admissible five-entry suffixes.
Using the GAP program ``automaton.g", which is available on arXiv.xxx, one obtains the eight continuation classes listed in
Table~\ref{tab:classes}. The program implements a partition refinement algorithm originally due to Edward F. Moore \cite{M56}. The program was in part constructed with the assistance of AI.
\end{proof}

\begin{table}[ht]
\centering
\caption{Continuation classes of the reachable five-entry suffixes.}
\label{tab:classes}
\renewcommand{\arraystretch}{1.25}
\begin{tabular}{c|l}
\hline
Class & Reachable five-entry suffixes \\
\hline
\(C_1\) &
\(
00000,\;
00100,\;
01000,\;
02000,\;
10000,\;
10100,\;
11000,\;
20000
\)
\\

\(C_2\) &
\(
00001,\;
00020,\;
01001,\;
10001,\;
10020,\;
11001,\;
20001
\)
\\

\(C_3\) &
\(
00002,\;
01002,\;
10002,\;
20002
\)
\\

\(C_4\) &
\(
00010,\;
00200,\;
01010,\;
10010,\;
11010,\;
20010
\)
\\

\(C_5\) &
\(
00011,\;
00201,\;
10011
\)
\\

\(C_6\) &
\(
00101,\;
01101,\;
02001,\;
10101,\;
20101
\)
\\

\(C_7\) &
\(
00110,\;
02010
\)
\\

\(C_8\) &
\(
01100,\;
20100
\)
\\
\hline
\end{tabular}
\end{table}

The construction depends only on the continuation classes not on any representatives.

For convenience we denote by
\[
q_1,\ldots,q_8
\]
the states corresponding to the continuation classes
\(C_1,\ldots,C_8\), respectively.

\subsection{The transition function}

The continuation classes introduced in the preceding subsection become
the states of a deterministic finite automaton only if appending a
frequency induces a well-defined transition between continuation
classes. We now verify that this is indeed the case.

\begin{proposition}\label{prop:welldefined}
Let \(u\) and \(v\) belong to the same continuation class.
For every symbol \(a\in\{0,1,2\}\), exactly one of the following two
possibilities occurs.

\begin{enumerate}
\item
Both words obtained by appending \(a\) to \(u\) and \(v\) are
inadmissible.

\item
Both words are admissible, and they belong to the same continuation
class.
\end{enumerate}

Consequently, appending a symbol induces a well-defined transition on
the set of continuation classes.
\end{proposition}

\begin{proof}
Suppose first that the extension \(ua\) is inadmissible. If \(va\)
were admissible, then the one-letter word \(a\) would be an admissible
continuation of \(v\), but not of \(u\), contradicting the assumption
that \(u\) and \(v\) belong to the same continuation class. Hence \(ua\)
and \(va\) are either both admissible or both inadmissible.

Now suppose that both extensions are admissible. Let \(w\) be an
arbitrary finite word over \(\Sigma\). Then
\[
uaw \text{ is admissible}
\quad\Longleftrightarrow\quad
vaw \text{ is admissible},
\]
because \(aw\) is a finite continuation and \(u\) and \(v\) belong to
the same continuation class. Since admissibility of future extensions depends only on the last five entries, the words obtained by appending \(a\) to \(u\) and \(v\) belong to the same continuation class.

 Hence the class reached after
appending \(a\) depends only on the continuation class of the original
suffix.
\end{proof}

Proposition~\ref{prop:welldefined} allows us to define a transition
function on the set of continuation classes.
Given a continuation class \(C_i\) and a symbol
\(a\in\{0,1,2\}\), we append \(a\) to any representative of \(C_i\).
If the resulting suffix is inadmissible, the transition leads to a
distinguished failure state.
Otherwise, the transition leads to the continuation class containing the
new suffix.
By Proposition~\ref{prop:welldefined}, this definition is independent of
the representative chosen.

We are therefore in a position to define the deterministic finite
automaton
\[
M_1=(Q_1,\Sigma,\delta,q_I,A),
\]
where
\[
Q_1=\{q_I,q_1,\ldots,q_8,q_f\},
\qquad
\Sigma=\{0,1,2\}.
\]
Here \(q_I\) is the initial state, \(q_i\) corresponds to the
continuation class \(C_i\), and \(q_f\) is the failure state. The set of
accepting states is
\[
A=\{q_I,q_1,\ldots,q_8\}.
\]
The transitions from the initial state are
\[
\delta(q_I,0)=q_1,
\qquad
\delta(q_I,1)=\delta(q_I,2)=q_f,
\]
and the remaining transitions are given in
Table~\ref{tab:transition}.

\begin{table}[ht]
\centering
\caption{Transition table of the deterministic finite automaton $M_1$.}
\label{tab:transition}
\renewcommand{\arraystretch}{1.2}
\begin{tabular}{c|ccc}
\hline
 & \(0\) & \(1\) & \(2\)\\
\hline
\(q_I\) & \(q_1\) & \(q_f\) & \(q_f\)\\
\hline
\(q_1\) & \(q_1\) & \(q_2\) & \(q_3\)\\
\(q_2\) & \(q_4\) & \(q_5\) & \(q_f\)\\
\(q_3\) & \(q_2\) & \(q_f\) & \(q_f\)\\
\(q_4\) & \(q_1\) & \(q_6\) & \(q_f\)\\
\(q_5\) & \(q_7\) & \(q_f\) & \(q_f\)\\
\(q_6\) & \(q_4\) & \(q_f\) & \(q_f\)\\
\(q_7\) & \(q_8\) & \(q_6\) & \(q_f\)\\
\(q_8\) & \(q_1\) & \(q_2\) & \(q_f\)\\
\(q_f\) & \(q_f\) & \(q_f\) & \(q_f\)\\
\hline
\end{tabular}
\end{table}

For prefixes of length less than five, the current suffix is understood
to be padded on the left by zeros, consistently with the initial window
\[
(0,0,0,0,0).
\]

\begin{theorem}\label{thm:automaton-correctness}
The deterministic finite automaton \(M_1\) accepts exactly the finite
frequency vectors satisfying conditions \textup{(F1)}--\textup{(F10)}.
Equivalently, it accepts exactly the frequency vectors of the partitions
satisfying conditions \textup{(IC1)} and \textup{(D1)}--\textup{(D10)}.
\end{theorem}

\begin{proof}
Let
\[
f_1f_2\cdots f_n
\]
be a finite word over \(\Sigma=\{0,1,2\}\).

Suppose first that the word is accepted by \(M_1\).
The initial transition forces \(f_1=0\), so condition \textup{(F1)}
holds. Since the alphabet is \(\Sigma\), condition \textup{(F2)} also
holds. At every subsequent step, the transition function sends the
current continuation class to the class determined by the new
five-entry suffix, unless appending the new symbol creates one of the
forbidden local configurations. In that case the transition enters the
failure state \(q_f\). Since the word is accepted, the computation never
enters \(q_f\). Hence none of the forbidden configurations in
Proposition~\ref{prop:forbidden} occurs, and the word is admissible.

Conversely, suppose that the word is admissible. Then \(f_1=0\), so the
first symbol sends the initial state \(q_I\) to \(q_1\). Assume
inductively that, after reading an admissible prefix, the automaton is in
the state corresponding to the continuation class of its current
five-entry suffix. Appending the next frequency cannot create a
forbidden local configuration, because the full word is admissible.
Therefore the transition does not enter \(q_f\); instead it leads to the
state corresponding to the continuation class of the new five-entry
suffix. By induction, the automaton never reaches \(q_f\), and after the
entire word has been read it lies in an accepting state. Thus the word
is accepted by \(M_1\).

The equivalence with conditions \textup{(IC1)} and \textup{(D1)}--\textup{(D10)} follows
from Proposition~\ref{prop:forbidden}.
\end{proof}

This completes the construction and verification of the finite
automaton associated with the admissible frequency vectors. 

\begin{remark}
By the Myhill-Nerode Theorem (see e.g. \cite{LP98} Theorem 2.5.2) there exists a unique, up to isomorphism, minimal automaton that accepts ${\bf f}(\mathcal{C})$. 

Using the algorithm in \cite{K97} Lecture 14 it is easy to check that $M_1$ is the minimal automaton that accepts ${\bf f}(\mathcal{C})$.
\end{remark}

\section[A linear system of q-difference equations]{A linear system of \texorpdfstring{$q$-}\ difference equations}\label{conj9}

\noindent
In this section, we associate a bivariate generating function to the set of partitions $\mathcal{C}$. Using the automaton $M_1$ from Section \ref{dfaconstruction} we determine a linear system of $q$-difference equations, and we find a solution to the system. This leads to a proof of the first combinatorial identity Theorem \ref{theorem1}.

For a partition $\lambda = (\lambda_1,\dots,\lambda_k)$, set $|\lambda| = \sum_{i=1}^k \lambda_i$, the {\em  weight} of $\lambda$, and $\ell(\lambda) = k$, the {\em length} of $\lambda$.
With $\mathcal{C}$  the set of partitions that satisfy the difference conditions (DC) and the initial condition (IC1), set
\begin{equation}
    h(x,q) = \sum_{\lambda \in\mathcal{C}}x^{\ell(\lambda)}q^{|\lambda|},
\end{equation}
which we consider as an element of $\mathbb{Z}[[x,q]]$. We include the empty partition, denoted $\epsilon$, as a partition in $\mathcal{C}$ with $\ell(\epsilon) = 0$ and $|\epsilon| = 0$, so that $h(x,q)\in 1+xq\mathbb{Z}[[x,q]]$. 

From now on we shall use the alternative notation $q_{11} = q_I$.
For $j=1,\dots,8,11$ let $\mathcal{C}_j$ be the set of partitions associated with frequency vectors that start at $q_j$ and have support in $A$. Set
\begin{equation}
    h_j(x,q)=\sum_{\lambda\in\mathcal{C}_j}x^{\ell(\lambda)}q^{|\lambda|}\in\mathbb{Z}[[x,q]].
\end{equation}

\begin{lemma}
    Assume that a state $q_a$ transitions to $q_b,q_c,q_d\in A$ by $0,1,2$ respectively.
    Then
    \[
    h_a(x,q) = h_b(xq,q) + xqh_c(xq,q) + x^2q^2h_d(xq,q).
    \]
    If instead $q_b,q_c\in A$ and $q_d=q_f$, then
    \[
    h_a(x,q) = h_b(xq,q) + xqh_c(xq,q).
    \]
    Finally, if $q_b\in A$ and $q_c=q_d=q_f$, then
    \[
    h_a(x,q) = h_b(xq,q).
    \]
\end{lemma}
\begin{proof}
    In the case $q_b,q_c,q_d\in A$, the set ${\bf f}(\mathcal{C}_a)$ decomposes as
    \begin{equation*}
        {\bf f}(\mathcal{C}_a) = \{(0,v)\mid v\in {\bf f}(\mathcal{C}_b)\}\cup \{(1,v)\mid v\in {\bf f}(\mathcal{C}_c)\}
        \cup \{(2,v)\mid v\in {\bf f}(\mathcal{C}_d)\}.
    \end{equation*}
Let $\lambda = (\lambda_1,\dots,\lambda_k)$ be a partition with frequency vector $v\in {\bf f}(\mathcal{C}_d)$. Then the partition with frequency vector $(2,v)$ is
\[
\mu = (\lambda_1+1,\dots,\lambda_k+1,1,1).
\]
This has $\ell(\mu) = \ell(\lambda)+2$ and $|\mu| = |\lambda| + \ell(\lambda) + 2$. Adding over all such $\mu$ gives
\[
\sum_{\mu} x^{\ell(\mu)}q^{|\mu|} = \sum_{\lambda\in\mathcal{C}_d} x^{\ell(\lambda) + 2}q^{|\lambda|+\ell(\lambda)+2} = x^2q^2h_d(xq,q).
\]
Similarly, partitions with frequency vectors $(1,v)$ with $v\in {\bf f}(\mathcal{C}_c)$ 
leads to the term $xqh_c(xq,q)$, and partitions with frequency vectors $(0,v)$ with  $v\in {\bf f}(\mathcal{C}_b)$ leads to the term $h_b(xq,q)$.

The proofs of the other two cases are analogous.
    
\end{proof}

By applying the lemma to the transitions in Table~\ref{tab:transition}, we obtain the equations 
\eqref{hh1}--\eqref{hh11} below:

\begin{equation}\label{hh1}
    h_1(x,q) = h_1(xq,q) + xqh_2(xq,q) + x^2q^2h_3(xq,q),
\end{equation}
\begin{equation}\label{hh2}
    h_2(x,q) = h_4(xq,q) + xqh_5(xq,q),
\end{equation}
\begin{equation}\label{hh3}
    h_3(x,q) = h_2(xq,q),
\end{equation}
\begin{equation}\label{hh4}
    h_4(x,q) = h_1(xq,q) + xqh_6(xq,q),
\end{equation}
\begin{equation}\label{hh5}
    h_5(x,q) = h_7(xq,q),
\end{equation}
\begin{equation}\label{hh6}
    h_6(x,q) = h_4(xq,q),
\end{equation}
\begin{equation}\label{hh7}
    h_7(x,q) = xqh_6(xq,q) + h_8(xq,q),
\end{equation}
\begin{equation}\label{hh8}
    h_8(x,q) = h_1(xq,q) + xqh_2(xq,q),
\end{equation}
\begin{equation}\label{hh11}
    h(x,q) = h_{11}(x,q) = h_1(xq,q).
\end{equation}
From \eqref{hh3},\eqref{hh5},\eqref{hh6},\eqref{hh8} we can solve for $h_3,h_5,h_6,h_8$ respectively in terms of $h_1,h_2,h_4,h_7$. Substituting for $h_3,h_5,h_6,h_8$ in \eqref{hh1},\eqref{hh2},\eqref{hh4},\eqref{hh7} leads to
\begin{equation}\label{e1}
    h_1(x,q) = h_1(xq,q) + xqh_2(xq,q) + x^2q^2h_2(xq^2,q),
\end{equation}
\begin{equation}\label{e2}
    h_2(x,q) = h_4(xq,q) + xqh_7(xq^2,q),
\end{equation}
\begin{equation}\label{e3}
    h_4(x,q) = h_1(xq,q) +  xqh_4(xq^2,q), 
\end{equation}
\begin{equation}\label{e4}
    h_7(x,q) = xqh_4(xq^2,q) + h_1(xq^2,q) + xq^2h_2(xq^2,q).
\end{equation}

We shall show that the system \eqref{e1}--\eqref{e4} is satisfied by
\begin{equation}\label{h1}
    h_1(x,q) = \sum_{i,j\ge 0}\frac{x^{i+2j}q^{i^2+3ij+4j^2}}{(q;q)_i(q;q)_j}(q^{-i-2j} - q^{-i-j} + q^{-j}),
\end{equation}
\begin{equation}\label{h2}
    h_2(x,q) = \sum_{i,j\ge 0}\frac{x^{i+2j}q^{i^2+3ij+4j^2}}{(q;q)_i(q;q)_j}q^{2j},
\end{equation}
\begin{equation}\label{h4}
    h_4(x,q) = \sum_{i,j\ge 0}\frac{x^{i+2j}q^{i^2+3ij+4j^2}}{(q;q)_i(q;q)_j},
\end{equation}
\begin{equation}\label{h7}
    h_7(x,q) = \sum_{i,j\ge 0}\frac{x^{i+2j}q^{i^2+3ij+4j^2}}{(q;q)_i(q;q)_j}q^{j}.
\end{equation}
Combining \eqref{h1} and \eqref{hh11} leads to
\begin{equation}\label{hxq}
    h(x,q) = \sum_{i,j\ge 0}\frac{x^{i+2j}q^{i^2+3ij+4j^2}}{(q;q)_i(q;q)_j}(1 - q^{j} + q^{i+j}).
\end{equation}
The double sum in \eqref{hxq} comes from \cite{AEH20} Theorem 3. The double sums \eqref{h2}--\eqref{h7} were obtained by analogy with \cite{AEH20} Remark 3.5.

For $a,b,c,d\in\mathbb{Z}$, define
\begin{equation}
    g_{a,b,c,d}(x,q) = x^d \sum_{i,j\ge 0}\frac{x^{i+2j}q^{i^2+3ij+4j^2}}{(q;q)_i(q;q)_j}q^{ai+bj+c}.
\end{equation}
(In most of the expressions $a,b,c,d$ will be nonnegative.)
The change of variables $x=t^2,i=k_2,j=k_1$ gives
\begin{equation*}
    g_{a,b,c,d}(x,q) = t^{2d}\sum_{k_1,k_2\ge 0}\frac{t^{4k_1+2k_2}q^{4k_1^2+3k_1k_2+k_2^2}}{(q;q)_{k_1}(q;q)_{k_2}} q^{ak_2+bk_1+c} = f_{b,a,c,2d}(t,q)
\end{equation*}
with $f_{a,b,c,d}$ the series D. Salazar defines in \cite{Sa24} page 408. The following is a special case of \cite{Sa24} Lemma 1.9 (iii), (iv).

\begin{lemma}\label{sa}
For $a,b,c,d\in\mathbb{Z}$,

\begin{equation*}
(i)\qquad g_{a,b,c,d}(x,q) - g_{a+1,b,c,d}(x,q) = g_{a+2,b+3,a+c+1,d+1}(x,q),
\end{equation*}

\begin{equation*}
(ii)\qquad g_{a,b,c,d}(x,q) - g_{a,b+1,c,d}(x,q) = g_{a+3,b+8,b+c+4,d+2}(x,q).
\end{equation*}
\end{lemma}

In Propositions~\ref{prop.e3}--\ref{prop.e4} we show that we have a solution to the system \eqref{e1}--\eqref{e4}.
\begin{proposition}\label{prop.e3}
The series $h_1,h_4$ as in \eqref{h1}, \eqref{h4} satisfies \eqref{e3}.
\end{proposition}
\begin{proof}
    By Lemma \ref{sa} 
\begin{equation*}
        xqh_4(xq^2,q) = xq\sum_{i,j\ge 0}\frac{x^{i+2j}q^{i^2+3ij+4j^2}}{(q;q)_i(q;q)_j}q^{2i+4j} 
        = g_{2,4,1,1} = g_{0,1,0,0} - g_{1,1,0,0}.
\end{equation*}
Hence
\begin{equation*}
        h_4(x,q) - xqh_4(xq^2,q) = \sum_{i,j\ge 0}\frac{x^{i+2j}q^{i^2+3ij+4j^2}}{(q;q)_i(q;q)_j}(1-q^j+q^{i+j}) =  h_1(xq,q).
\end{equation*}
\end{proof}

\begin{proposition}\label{prop.e1}
The series $h_1,h_2$ as in \eqref{h1}, \eqref{h2} satisfies \eqref{e1}.
\end{proposition}
\begin{proof}
We shall show that
\begin{equation*}
    -h_1(x,q) + h_1(xq,q) + xqh_2(xq,q) + x^2q^2h_2(xq^2,q) = 0.
\end{equation*}
The first two terms transform as
\begin{equation*}
    \begin{split}
        &-h_1(x,q) + h_1(xq,q) \\
        &= -\sum_{i,j\ge 0}\frac{x^{i+2j}q^{i^2+3ij+4j^2}}{(q;q)_i(q;q)_j}
        (q^{-i-2j} - q^{-i-j} + q^{-j}) + \sum_{i,j\ge 0}\frac{x^{i+2j}q^{i^2+3ij+4j^2}}{(q;q)_i(q;q)_j}(1-q^j + q^{i+j}) \\
        &= (g_{0,0,0,0}-g_{-1,0,0,0}) + (g_{-1,0,0,0}-g_{-1,-2,0,0}) - (g_{0,1,0,0}-g_{1,1,0,0}) + (g_{-1,-1,0,0}-g_{0,-1,0,0}) \\
        &= -g_{1,3,0,1} + (g_{-1,0,0,0}-g_{-1,-1,0,0}) + (g_{-1,-1,0,0}-g_{-1,-2,0,0})
        -g_{2,4,1,1} + g_{1,2,0,1} \\
        &= (g_{1,2,0,1}-g_{1,3,0,1}) - g_{2,7,3,2} - g_{2,6,2,2} - g_{2,4,1,1} \\
        &= g_{4,10,6,3} - g_{2,7,3,2} - g_{2,6,2,2} - g_{2,4,1,1}.
    \end{split}   
\end{equation*}
Furthermore,
\begin{equation*}
    xqh_2(xq,q) = xq\sum_{i,j\ge 0}\frac{x^{i+2j}q^{i^2+3ij+4j^2}}{(q;q)_i(q;q)_j}q^{i+4j} = g_{1,4,1,1},
\end{equation*}
and
\begin{equation*}
    x^2q^2h_2(xq^2,q) = x^2q^2\sum_{i,j\ge 0}\frac{x^{i+2j}q^{i^2+3ij+4j^2}}{(q;q)_i(q;q)_j}q^{2i+6j} = g_{2,6,2,2}.
\end{equation*}
Hence
\begin{equation*}
    \begin{split}
        &-h_1(x,q) + h_1(xq,q) + xqh_2(xq,q) + x^2q^2h_2(xq^2,q) \\
        &= g_{4,10,6,3} - g_{2,7,3,2} + (g_{1,4,1,1}-g_{2,4,1,1}) \\
        &=g_{4,10,6,3} - (g_{2,7,3,2} - g_{3,7,3,2}) = g_{4,10,6,3} - g_{4,10,6,3} = 0.
    \end{split}
\end{equation*}    
\end{proof}

\begin{proposition}\label{prop.e2}
The series $h_2,h_4,h_7$ as in \eqref{h2}, \eqref{h4}, \eqref{h7} satisfies \eqref{e2}.
\end{proposition}
\begin{proof}
By Lemma \ref{sa},
\begin{equation*}
\begin{split}
    h_2(x,q) - h_4(xq,q) &= \sum_{i,j\ge 0}\frac{x^{i+2j}q^{i^2+3ij+4j^2}}{(q;q)_i(q;q)_j}(q^{2j} - q^{i+2j}) 
    = g_{0,2,0,0} - g_{1,2,0,0} = g_{2,5,1,1}.
    \end{split}
\end{equation*}
    On the other hand
\begin{equation*}
        xqh_7(xq^2,q) = xq\sum_{i,j\ge 0}\frac{x^{i+2j}q^{i^2+3ij+4j^2}}{(q;q)_i(q;q)_j}q^{2i+5j}
        = g_{2,5,1,1}.
\end{equation*}
\end{proof}

\begin{proposition}\label{prop.e4}
The series $h_1,h_2,h_4,h_7$ as in \eqref{h1}--\eqref{h7} satisfies \eqref{e4}.
\end{proposition}
\begin{proof}
We shall show that
\begin{equation*}
    -h_7(x,q) + xqh_4(xq^2,q) + h_1(xq^2,q) + xq^2h_2(xq^2,q) = 0.
\end{equation*}
We obtain
\begin{equation*}
	\begin{split}
		&-h_7(x,q) + xqh_4(xq^2,q) + h_1(xq^2,q) + xq^2h_2(xq^2,q) \\
		&= -\sum_{i,j\ge 0}\frac{x^{i+2j}q^{i^2 + 3ij + 4j^2}}{(q;q)_i(q;q)_j}q^j
		 + x\sum_{i,j\ge 0}\frac{x^{i+2j}q^{i^2 + 3ij + 4j^2}}{(q;q)_i(q;q)_j}q^{2i+4j+1} \\
			&\quad + \sum_{i,j\ge 0}\frac{x^{i+2j}q^{i^2 + 3ij + 4j^2}}{(q;q)_i(q;q)_j}q^{i+2j}(1-q^j+q^{i+j})
		 + x\sum_{i,j\ge 0}\frac{x^{i+2j}q^{i^2 + 3ij + 4j^2}}{(q;q)_i(q;q)_j}q^{2i+6j+2} \\
		&=  -(g_{0,1,0,0}-g_{1,1,0,0}) - (g_{1,1,0,0}-g_{1,2,0,0})-(g_{1,3,0,0}-g_{2,3,0,0})+g_{2,4,1,1} + g_{2,6,2,1}    \\
		&=  -g_{2,4,1,1} -g_{4,9,5,2} + (g_{2,6,2,1}- g_{3,6,2,1}) + g_{2,4,1,1}    \\
		&= -g_{4,9,5,2} + g_{4,9,5,2} = 0.   
	\end{split}
\end{equation*}
\end{proof}

As $h_1(x,q)\in 1+xq\mathbb{Z}[[x,q]]$, we have $h_1(0,q) = 1$.
We shall show that the system \eqref{e1}--\eqref{e4} has a unique solution satisfying the initial condition $h_1(0,q) = 1$. This is analogous to \cite{AEH20} Lemma 3.2.
Introduce the coefficients of the formal power series $h_1(x,q)$, $h_2(x,q)$, $h_4(x,q)$, $h_7(x,q)$ by
	\begin{equation*}
		h_1(x,q) = \sum_{m,n\ge 0} a(m,n)x^mq^n, 
	\end{equation*}
	\begin{equation*}
		h_2(x,q) = \sum_{m,n\ge 0} b(m,n)x^mq^n, 
	\end{equation*}
	\begin{equation*}
		h_4(x,q) = \sum_{m,n\ge 0} c(m,n)x^mq^n,
	\end{equation*}
	\begin{equation*}
		h_7(x,q) = \sum_{m,n\ge 0} d(m,n)x^mq^n.
	\end{equation*}
	Then the system \eqref{e1}--\eqref{e4} is equivalent to
	\begin{eqnarray}
		a(m,n)&=& a(m,n-m) + b(m-1,n-m) + b(m-2,n-2m+2), \label{eq9} \\
		b(m,n)&=& c(m,n-m) + d(m-1,n-2m+1), \label{eq10}\\
		c(m,n)&=& a(m,n-m) + c(m-1,n-2m+1), \label{eq11}\\
		d(m,n)&=& c(m-1,n-2m+1) + a(m,n-2m) + b(m-1,n-2m).\label{eq12}
	\end{eqnarray}
	Here, for $x\in\{a,b,c,d\}$, we have extended the notation $x(m,n)$ to all $(m,n)\in\mathbb{Z}^2$ by setting $x(m,n) = 0$ if $m<0$ or $n<0$. The initial condition is
	\begin{align*}
		a(0,0) &= 1, \\
		a(0,n) &= 0 \text{\ for\ } n\ge 1.
	\end{align*} 
	Substituting $m=0$ in \eqref{eq9}--\eqref{eq12} gives $a(0,n) = b(0,n) = c(0,n) = d(0,n)$, so it follows that
	\begin{align*}
		b(0,0) &= c(0,0) = d(0,0) = 1, \\
		b(0,n) &= c(0,n) = d(0,n) = 0 \text{\ for\ } n\ge 1.
	\end{align*}

\begin{proposition}\label{prop.unique}
		The series $h_1,h_2,h_4,h_7$ are uniquely determined by the system \eqref{e1}--\eqref{e4} together with the initial condition $h_1(0,q) = 1$.
	\end{proposition}
	\begin{proof}
		We shall show that $a(m,n), b(m,n), c(m,n), d(m,n)$ are uniquely determined by \eqref{eq9}--\eqref{eq12} together with the above initial condition. By the above arguments, we may restrict ourselves to $m\ge 1$. We prove the uniqueness by induction on $n$. Let $S_n$ be the assertion:
	\[
	\{a(m,k),b(m,k),c(m,k),d(m,k)\}
	\]
	are uniquely determined for all pairs $(m,k)$ with $m\in\mathbb{Z}$, $k\le n$.
	To start the induction, note that $S_{-1}$ holds. Assume that $S_{n-1}$ is true for some $n\ge 0$. We show that $S_n$ holds. By \eqref{eq10}, \eqref{eq11}, \eqref{eq12} $b(m,n), c(m,n), d(m,n)$ are unique. If $m=1$, then  
    $b(m-2,n-2m+2)=b(-1,n)=0$, so that uniqueness of $a(m,n)$ follows from \eqref{eq9}. If $m\ge 2$ uniqueness of $a(m,n)$ follows from \eqref{eq9} and the induction hypothesis.
	\end{proof}

Combining \eqref{hxq} with results from \cite{AEH20}, we obtain the following
\begin{theorem}\label{theorem1}
The generating function of the partitions in $\mathcal{C}$ satisfies
\begin{equation}\label{nahmseries}
    \sum_{\lambda\in\mathcal{C}} x^{\ell(\lambda)}q^{|\lambda|} = \sum_{i,j\ge 0}\frac{x^{i+2j}q^{i^2+3ij+4j^2}}{(q;q)_i(q;q)_j}(1-q^j + q^{i+j}).
\end{equation}
In particular,
\begin{equation}
    \sum_{\lambda\in\mathcal{C}} q^{|\lambda|} = \sum_{i,j\ge 0}\frac{q^{i^2+3ij+4j^2}}{(q;q)_i(q;q)_j}(1-q^j + q^{i+j})
= \prod_{n\equiv \pm2,\pm3,\pm4,\pm5\pmod{16}}(1-q^n)^{-1}.
\end{equation}
\end{theorem} 
\begin{proof}
    By Propositions \ref{prop.e3} -- \ref{prop.unique} the equations \eqref{h1}--\eqref{h7} hold for the original bi-variate generating functions. Hence \eqref{hxq} holds, and this gives \eqref{nahmseries}. The identity
    \begin{equation*}
    \sum_{i,j\ge 0}\frac{q^{i^2+3ij+4j^2}}{(q;q)_i(q;q)_j}(1-q^j + q^{i+j})
= \prod_{n\equiv \pm2,\pm3,\pm4,\pm5\pmod{16}}(1-q^n)^{-1}
\end{equation*} 
is proved in \cite{AEH20} Theorem 3.
\end{proof}

\section[The set of partitions D]{The set of partitions \texorpdfstring{$\mathcal{D}$}{D}}
\noindent
Let $\mathcal{D}$ denote the set of integer partitions $\lambda$ such that the difference conditions (DC) and the initial conditions (IC2) hold, where (IC2) are the conditions that $\lambda$ does not have a sub-partition
$(1,1)$, or $(2,2)$. 

Let $(f_1,f_2,\dots)$ be the frequency vector associated with $\lambda$ as in 
Section~\ref{frequency_vectors}. Then conditions (DC) and (IC2) translate as follows to conditions on $(f_j)_{j=1}^\infty$:

\begin{description}

    \item[1] 
    conditions (F2)--(F10) hold,
    
    \item[2]
    $f_1\in \{0,1\}$, 
    
    \item[3]
    $f_2\in \{0,1\}$. 
    
\end{description}

The set of frequency vectors associated to $\mathcal{D}$ is accepted by a DFA $M_2 = (Q_2, \Sigma, \delta, s, A)$ where the set of states is $Q_2 = \{q_1,\dots,q_8,q_f,q_{12}\}$, the set of transitions is 
$\Sigma = \{0,1,2\}$, the transition function $\delta$ is given in Table~\ref{tab:tb3}, the start state is $s=q_{12}$, and the set of ``accept states" is $A = Q_2\setminus \{q_f\}$. The set of finite paths starting at $s$ and with support in $A$ corresponds to the partitions in $\mathcal{D}$.

\begin{table}[ht]
\centering
\caption{Transition table of the deterministic finite automaton $M_2$.}
\label{tab:tb3}
\renewcommand{\arraystretch}{1.2}
\begin{tabular}{c|ccc}
\hline
 & \(0\) & \(1\) & \(2\)\\
\hline
\(q_{12}\) & \(q_8\) & \(q_2\) & \(q_f\)\\
\hline
\(q_1\) & \(q_1\) & \(q_2\) & \(q_3\)\\
\(q_2\) & \(q_4\) & \(q_5\) & \(q_f\)\\
\(q_3\) & \(q_2\) & \(q_f\) & \(q_f\)\\
\(q_4\) & \(q_1\) & \(q_6\) & \(q_f\)\\
\(q_5\) & \(q_7\) & \(q_f\) & \(q_f\)\\
\(q_6\) & \(q_4\) & \(q_f\) & \(q_f\)\\
\(q_7\) & \(q_8\) & \(q_6\) & \(q_f\)\\
\(q_8\) & \(q_1\) & \(q_2\) & \(q_f\)\\
\(q_f\) & \(q_f\) & \(q_f\) & \(q_f\)\\
\hline
\end{tabular}
\end{table}

The automaton $M_2$ can be obtained by the GAP program with the Automata package as in \cite{TT19} Remark 4.3.
Compared to Section \ref{dfaconstruction}, the state $q_{11}$ has been replaced by a state $q_{12}$. For $j=1,\dots,8,12$, denote by $\mathcal{D}_j$ the set of partitions corresponding to frequency vectors that start at $q_j$ and have support in $A$. Set $h_j(x,q) = \sum_{\lambda\in\mathcal{D}_j}x^{\ell(\lambda)}q^{|\lambda|}$. Then, for $j=1,\dots,8$, we have  $\mathcal{C}_j = \mathcal{D}_j$, so that, for $j=1,\dots,8$, the $h_j$ are given by the expressions in Section \ref{conj9}.

\begin{proposition}\label{bivarD}
    The series $h_{12}$ is given by
    \begin{equation}
         h_{12}(x,q) = \sum_{i,j\ge 0}\frac{x^{i+2j}q^{i^2+3ij+4j^2}}{(q;q)_i(q;q)_j}q^{i+j} + x\sum_{i,j\ge 0}\frac{x^{i+2j}q^{i^2+3ij+4j^2}}{(q;q)_i(q;q)_j}q^{i+4j+1}.
    \end{equation}
\end{proposition}
\begin{proof}
    From the automaton, we have 
    \begin{equation*}
        h_{12}(x,q) = h_8(xq,q) + xqh_2(xq,q).
    \end{equation*}
By \eqref{hh8},
\begin{equation*}
    h_8(xq,q) = h_1(xq^2,q) + xq^2h_2(xq^2,q),
\end{equation*}
so that
\begin{equation*}
    h_{12}(x,q) = h_1(xq^2,q) + xqh_2(xq,q) + xq^2h_2(xq^2,q).
\end{equation*}
The expressions \eqref{h1}, \eqref{h2} give
\begin{equation*}
    \begin{split}
        h_{12}(x,q) &= \sum_{i,j\ge 0}\frac{x^{i+2j}q^{i^2+3ij+4j^2}}{(q;q)_i(q;q)_j}q^{i+2j}(1-q^j+q^{i+j})
        +x\sum_{i,j\ge 0}\frac{x^{i+2j}q^{i^2+3ij+4j^2}}{(q;q)_i(q;q)_j}q^{i+4j+1} \\
        &\quad +x\sum_{i,j\ge 0}\frac{x^{i+2j}q^{i^2+3ij+4j^2}}{(q;q)_i(q;q)_j}q^{2i+6j+2} \\
        &= x\sum_{i,j\ge 0}\frac{x^{i+2j}q^{i^2+3ij+4j^2}}{(q;q)_i(q;q)_j}q^{i+4j+1} 
        +g_{1,2,0,0} - g_{1,3,0,0} + g_{2,3,0,0} + g_{2,6,2,1}.
    \end{split}
\end{equation*}
With Lemma~\ref{sa} the sum of the last four terms transforms as follows:
\begin{equation*}
    \begin{split}
    &g_{1,2,0,0} - (g_{1,3,0,0} - g_{2,3,0,0}) + g_{2,6,2,1} 
    = g_{1,1,0,0}- (g_{1,1,0,0}-g_{1,2,0,0}) + (g_{2,6,2,1}- g_{3,6,2,1}) \\
    &= g_{1,1,0,0} - g_{4,9,5,2} + g_{4,9,5,2} = g_{1,1,0,0}.
    \end{split}
\end{equation*}
\end{proof}

Combining Proposition~\ref{bivarD} with identities from \cite{AEH20} we obtain the following
\begin{theorem}\label{theorem2}
    The generating function of the partitions in $\mathcal{D}$ satisfies
    \begin{equation}
        \sum_{\lambda\in\mathcal{D}}x^{\ell(\lambda)}q^{|\lambda|}= \sum_{i,j\ge 0}\frac{x^{i+2j}q^{i^2+3ij+4j^2}}{(q;q)_i(q;q)_j}q^{i+j} + x\sum_{i,j\ge 0}\frac{x^{i+2j}q^{i^2+3ij+4j^2}}{(q;q)_i(q;q)_j}q^{i+4j+1}.
    \end{equation}
    In particular,
    \begin{equation}
        \sum_{\lambda\in\mathcal{D}}q^{|\lambda|} = \sum_{i,j\ge 0}\frac{q^{i^2+3ij+4j^2}}{(q;q)_i(q;q)_j}
        (q^{i+j} + q^{i+4j+1}) = \prod_{n=1}^\infty (1+q^n).
    \end{equation}
\end{theorem}
\begin{proof}
    As $q_{12}$ is the start state of the automaton $M_2$, $h_{12}$ gives the bi-variate generating function of the partitions in $\mathcal{D}$. The identity
    \begin{equation}
        \sum_{i,j\ge 0}\frac{q^{i^2+3ij+4j^2}}{(q;q)_i(q;q)_j}
        (q^{i+j} + q^{i+4j+1}) = \prod_{n=1}^\infty (1+q^n)
    \end{equation}
    results from the two expressions (1.12.1) and (1.12.2) in \cite{AEH20} for the character of the irreducible highest weight module of highest weight $(c,h) = (1/2,1/16)$ for the Virasoro Lie algebra.
\end{proof}

\section[The set of partitions E]{The set of partitions \texorpdfstring{$\mathcal{E}$}{E}}
\noindent
Let $\mathcal{E}$ denote the set of integer partitions $\lambda$ such that the difference conditions (DC) and the initial conditions (IC3) hold, where (IC3) are the conditions that $\lambda$ does not have a sub-partition
$(1,1)$, $(2,1)$, $(2,2)$, $(3,2)$, $(3,3)$, $(4,4,1)$. 

Let $(f_1,f_2,\dots)$ be the frequency vector associated with $\lambda$ as in Section~\ref{frequency_vectors}. In Proposition \ref{prop:IC3} it will be shown that conditions (DC) and (IC3) translate as follows to conditions on $(f_j)_{j=1}^\infty$:

\begin{description}

    \item[1] 
     conditions (F2)--(F10) hold,
    
    \item[2]
    $f_1\in \{0,1\}$, 
    
    \item[3]
    $(f_1,f_2) \neq (1,1)$,

    \item[4]
    $(f_1,f_2) \neq (0,2)$,

    \item[5]
    $(f_1,f_2,f_3) \neq  (0,1,1)$,

    \item[6]
    $(f_1,f_2,f_3) \neq  (0,0,2)$,

    \item[7]
    $(f_1,f_2,f_3,f_4) \neq  (1,0,0,2)$.
    
\end{description}

The set of frequency vectors associated to $\mathcal{E}$ is accepted by a DFA $M_3 = (Q_3, \Sigma, \delta, s, A)$ where the set of states is $Q_3 = \{q_1,\dots,q_8,q_f,q_{13}\}$, the set of transitions is $\Sigma = \{0,1,2\}$, the transition function $\delta$ is given in Table~\ref{tab:tb4}, the start state is $s=q_{13}$, and the set of ``accept states" is $A = Q_3\setminus \{q_f\}$. The set of finite paths starting at $s$ and with support in $A$ corresponds to the partitions in $\mathcal{E}$.

\begin{table}[ht]
\centering
\caption{Transition table of the deterministic finite automaton $M_3$.}
\label{tab:tb4}
\renewcommand{\arraystretch}{1.2}
\begin{tabular}{c|ccc}
\hline
 & \(0\) & \(1\) & \(2\)\\
\hline
\(q_{13}\) & \(q_7\) & \(q_5\) & \(q_f\)\\
\hline
\(q_1\) & \(q_1\) & \(q_2\) & \(q_3\)\\
\(q_2\) & \(q_4\) & \(q_5\) & \(q_f\)\\
\(q_3\) & \(q_2\) & \(q_f\) & \(q_f\)\\
\(q_4\) & \(q_1\) & \(q_6\) & \(q_f\)\\
\(q_5\) & \(q_7\) & \(q_f\) & \(q_f\)\\
\(q_6\) & \(q_4\) & \(q_f\) & \(q_f\)\\
\(q_7\) & \(q_8\) & \(q_6\) & \(q_f\)\\
\(q_8\) & \(q_1\) & \(q_2\) & \(q_f\)\\
\(q_f\) & \(q_f\) & \(q_f\) & \(q_f\)\\
\hline
\end{tabular}
\end{table}

The automaton $M_3$ can be obtained by the GAP program with the Automata package as in \cite{TT19} Remark 4.3.
Compared to Section~\ref{dfaconstruction}, the state $q_{11}$ has been replaced by a state $q_{13}$. For $j=1,\dots,8,13$, denote by $\mathcal{E}_j$ the set of partitions corresponding to frequency vectors that start at $q_j$ and have support in $A$. Set $h_j(x,q) = \sum_{\lambda\in\mathcal{E}_j}x^{\ell(\lambda)}q^{|\lambda|}$. Then, for $j=1,\dots,8$, we have  $\mathcal{C}_j = \mathcal{E}_j$, so that, for $j=1,\dots,8$, the $h_j$ are given by the expressions in Section \ref{conj9}.
\begin{proposition}\label{bivarE}
    The series $h_{13}$ is given by
    \begin{equation}\label{h13}
    \begin{split}
         h_{13}(x,q) &= \sum_{i,j\ge 0}\frac{x^{i+2j}q^{i^2+3ij+4j^2}}{(q;q)_i(q;q)_j}q^{i+3j} + x\sum_{i,j\ge 0}\frac{x^{i+2j}q^{i^2+3ij+4j^2}}{(q;q)_i(q;q)_j}q^{2i+5j+1}.
    \end{split}
    \end{equation}
\end{proposition}
\begin{proof}
From the automaton and \eqref{hh5} we have
\begin{equation*}
    \begin{split}
        h_{13}(x,q) &= h_7(xq,q) + xqh_5(xq,q) = h_7(xq,q) + xqh_7(xq^2,q) \\
        &= \sum_{i,j\ge 0}\frac{x^{i+2j}q^{i^2+3ij+4j^2}}{(q;q)_i(q;q)_j}q^{i+3j} + x\sum_{i,j\ge 0}\frac{x^{i+2j}q^{i^2+3ij+4j^2}}{(q;q)_i(q;q)_j}q^{2i+5j+1}.
    \end{split}
\end{equation*}
\end{proof}

Define $q$-series $g_{a,b,c}$, obtained by substituting $x=1$ in $g_{a,b,c,d}(x,q)$:
\begin{equation*}
    g_{a,b,c}(q) = \sum_{i,j \ge 0}\frac{q^{i^2 + 3ij +4j^2}}{(q;q)_i(q;q)_j}q^{ai+bj+c}.
\end{equation*}
Substituting $x=1$ in Lemma~\ref{sa} gives
\begin{lemma}\label{sa2}
    For $a,b,c\in\mathbb{Z}$,

\begin{equation*}
(i)\qquad g_{a,b,c}(q) - g_{a+1,b,c}(q) = g_{a+2,b+3,a+c+1}(q),
\end{equation*}

\begin{equation*}
(ii)\qquad g_{a,b,c}(q) - g_{a,b+1,c}(q) = g_{a+3,b+8,b+c+4}(q).
\end{equation*}
\end{lemma}

The formulas (1.12.1) and (1.12.2) in \cite{AEH20} for the character of the irreducible Virasoro module $V_{1/2}$ with the highest weight $(c,h) = (1/2,1/2)$ give
\begin{equation*}
    \begin{split}
        \frac{1}{2}&\left(\prod_{n=1}^\infty (1+q^{n-1/2}) - \prod_{n=1}^\infty (1-q^{n-1/2}) \right)
= \chi_{V_{1/2}}(q) = q^{1/2} \sum_{i,j \ge 0}\frac{q^{i^2 + 3ij +4j^2}}{(q;q)_i(q;q)_j}(q^{2j} - q^{4i + 10j + 6}) 
\\
&= q^{1/2} (g_{0,2,0} - g_{4,10,6}).
    \end{split}
\end{equation*}
The case $(a,b,c) = (1,2,0)$ of Lemma~\ref{sa2}(ii) gives
\[
g_{1,2,0} - g_{1,3,0} = g_{4,10,6}
\]
so that
\[
g_{0,2,0} - g_{4,10,6} = g_{0,2,0} - g_{1,2,0} + g_{1,3,0}.
\]
In addition,
\[
g_{0,2,0}-g_{1,2,0} = g_{2,5,1}
\]
so that
\[
g_{0,2,0} - g_{4,10,6} = g_{2,5,1} + g_{1,3,0}.
\]
Thus, substituting $x=1$ in \eqref{h13} gives
\begin{equation}\label{firsteexpression}
    \begin{split}
\sum_{\mu\in\mathcal{E}}q^{|\mu|} &= \sum_{i,j \ge 0}\frac{q^{i^2 + 3ij +4j^2}}{(q;q)_i(q;q)_j}(q^{i+3j} + q^{2i+5j+1}) = g_{1,3,0} + g_{2,5,1} = g_{0,2,0} - g_{4,10,6} \\
&=\frac{1}{2}q^{-1/2}\left(\prod_{n=1}^\infty (1+q^{n-1/2}) - \prod_{n=1}^\infty (1-q^{n-1/2}) \right).
\end{split}
\end{equation}

\begin{proposition}\label{prop6.3}
    We have
    \begin{equation}\label{modularf}
    \frac{1}{2}q^{-1/2}\left(\prod_{n=1}^\infty (1+q^{n-1/2}) 
    - \prod_{n=1}^\infty (1-q^{n-1/2})\right)
    = \prod_{n\equiv \pm1,\pm4,\pm6,\pm7 \pmod{16}}(1-q^n)^{-1}.
    \end{equation}
\end{proposition}
\begin{proof}
We shall use the Jacobi triple product identity \cite{A76} Theorem 2.8
\begin{equation}\label{jacobitriple}
\prod_{n=1}^\infty \{(1+xq^{2n-1})(1+x^{-1}q^{2n-1})(1-q^{2n})\} = \sum_{k\in\mathbb{Z}}x^kq^{k^2}.
\end{equation}
Substitute first $x=q^{1/2}$ to obtain
\begin{equation}\label{plusproduct}
\prod_{n=1}^\infty \{(1+q^{n-1/2})(1-q^{2n})\} = \sum_{k\in\mathbb{Z}}q^{k^2+k/2}.
\end{equation}
Similarly, substituting $x=-q^{1/2}$ in \eqref{jacobitriple} gives
\begin{equation}\label{minusproduct}
\prod_{n=1}^\infty \{(1-q^{n-1/2})(1-q^{2n})\} = \sum_{k\in\mathbb{Z}}(-1)^kq^{k^2+k/2}.
\end{equation}
Subtract \eqref{minusproduct} from \eqref{plusproduct} to obtain
\begin{equation}\label{lhstransf}
    \begin{split}
&\prod_{n=1}^\infty (1-q^{2n})\left(\prod_{n=1}^\infty (1+q^{n-1/2}) - \prod_{n=1}^\infty (1-q^{n-1/2})  \right) = 2\sum_{k\in\mathbb{Z}}q^{(2k-1)^2 + (2k-1)/2} \\
&= 2q^{1/2}\sum_{k\in\mathbb{Z}} q^{4k^2-3k}.
\end{split}
\end{equation}
It is easy to transform
\begin{equation}\label{producttransf}
\begin{split}
&\prod_{n=1}^\infty (1-q^{2n}) \prod_{n\equiv \pm1,\pm4,\pm6,\pm7 \pmod{16}}(1-q^n)^{-1} 
= \prod_{n=1}^\infty \{(1+q^{8n-7})(1+q^{8n-1})(1-q^{8n})\} \\
&= \sum_{k\in\mathbb{Z}} q^{4k^2-3k}.
\end{split}
\end{equation}
In the last equality, we applied \eqref{jacobitriple} with $q$ replaced by $q^4$ and $x$ replaced by $q^{-3}$. Now \eqref{modularf} follows from \eqref{lhstransf} and \eqref{producttransf}.
\end{proof}
Proposition~\ref{prop6.3} can alternatively be proved by modifying the argument in \cite{B51} p.220.
Combining \eqref{firsteexpression} with Proposition~\ref{prop6.3}, we have proved the following result: 

\begin{theorem}\label{theorem3}
    The generating function of the partitions in $\mathcal{E}$ satisfies
    \begin{equation}
        \sum_{\lambda\in\mathcal{E}}x^{\ell(\lambda)}q^{|\lambda|}= \sum_{i,j\ge 0}\frac{x^{i+2j}q^{i^2+3ij+4j^2}}{(q;q)_i(q;q)_j}q^{i+3j} + x\sum_{i,j\ge 0}\frac{x^{i+2j}q^{i^2+3ij+4j^2}}{(q;q)_i(q;q)_j}q^{2i+5j+1}.
    \end{equation}
    In particular,
    \begin{equation}
        \begin{split}
       & \sum_{\lambda\in\mathcal{E}}q^{|\lambda|} = \sum_{i,j\ge 0}\frac{q^{i^2+3ij+4j^2}}{(q;q)_i(q;q)_j}
        (q^{i+3j} + q^{2i+5j+1}) \\
        &= \frac{1}{2}q^{-1/2}\left(\prod_{n=1}^\infty (1+q^{n-1/2}) - \prod_{n=1}^\infty (1-q^{n-1/2}) \right) = \prod_{n\equiv \pm1,\pm4,\pm6,\pm7 \pmod{16}}(1-q^n)^{-1}.
        \end{split}
    \end{equation}
\end{theorem}

\section{Appendix}
\label{sec:appendix}

\subsection*{Translation of the (IC3) conditions}

In this Appendix we reformulate the (IC3)  partition conditions as local
conditions on the frequency vector.

\begin{proposition}\label{prop:IC3}
  A partition satisfies conditions (IC3) and (D1)--(D10) if and only if its frequency vector satisfies conditions (F2)--(F10) together with
$$
  \begin{matrix}
    f_1\in \{0,1\} \\
    (f_1,f_2)\neq (1,1) \\
    (f_1,f_2)\neq (0,2) \\
    (f_1,f_2, f_3)\neq (0,1,1)\\
    (f_1,f_2, f_3)\neq (0,0,2) \\
    (f_1,f_2, f_3,f_4)\neq (1,0,0,2)
  \end{matrix}
$$

\end{proposition}

\begin{proof}
We already showed that the difference conditions (D1)--(D10) are equivalent to the frequency conditions (F2)--(F10). We only need to translate the six forbidden subpartitions in (IC3).

{\bf Step 1}
The subpartition $(1,1)$ occurs if and only if $f_1\geq 2$, but since $f_1\in \{0,1,2\}$, this is equivalent to $f_1=2$.  Thus excluding the subpartition $(1,1)$ is equivalent to requiring $f_1\in\{0,1\}$.

{\bf Step 2}
The subpartition $(2,1)$ occurs if and only if $f_1\geq 1, f_2\geq 1$, so the possible configurations are
$$
(f_1,f_2)=\begin{cases}
  (1,1)\\
  (1,2)\\
  (2,1)\\
  (2,2)
\end{cases}
$$
The last three are already forbidden by (F3). Hence the only new frequency configuration that must be excluded is $(f_1,f_2)=(1,1)$.

{\bf Step 3}
The subpartition $(2,2)$ occurs if and only if $ f_2\geq 2$, which forces $f_2=2$, so the possible configurations are
$$
(f_1,f_2)=\begin{cases}
  (0,2)\\
  (1,2)\\
    (2,2)
\end{cases}
$$
The last two are  forbidden by (F3). Hence the only new frequency configuration that must be excluded is $(f_1,f_2)=(0,2)$.

{\bf Step 4}
The subpartition $(3,2)$ occurs if and only if $ f_2\geq 1, f_3\geq 1$ so the possible configurations for the pair $(f_2,f_3)$ are
$$
(f_2,f_3)=\begin{cases}
  (1,1)\\
  (1,2)\\
  (2,1)\\
    (2,2)
\end{cases}
$$
The last three are  forbidden by (F3). Hence the only new frequency configuration that must be excluded is $(f_2,f_3)=(1,1)$. Moreover, since $f_1\in\{0,1,2\}$,
the remaining possibilities are 
$$
(f_1,f_2,f_3)=\begin{cases}
  (0,1,1)\\
  (1,1,1)\\
  (2,1,1)
\end{cases}
$$
But $(f_1,f_2,f_3)=(1,1,1)$ is already forbidden by (F4) and $(f_1,f_2,f_3)=(2,1,1)$
 is forbidden by (F3) (in the first two entries). Hence the only new frequency configuration to be excluded is $(0,1,1)$.

{\bf Step 5}
The subpartition $(3,3)$ occurs if and only if $ f_3\geq 2$ which means $f_3=2$. The possible configurations for the pair $(f_2,f_3)$ are
$$
(f_2,f_3)=\begin{cases}
  (0,2)\\
  (1,2)\\
     (2,2)
\end{cases}
$$
The last two are  forbidden by (F3), so $f_2=0$. Next, $f_1\in \{0,1,2\}$ so
the remaining possibilities are
$$
(f_1,f_2,f_3)=\begin{cases}
  (0,0,2)\\
  (1,0,2)\\
  (2,0,2)
\end{cases}
$$
But  (F5) forbids the last two. Thus the only new frequency configuration to be excluded is $(0,0,2)$.

{\bf Step 6}
Finally the subpartition $(4,4,1)$ occurs if and only if $f_1\geq 1, f_4\geq 2$ if and only if $f_1\in \{1,2\}, f_4=2$. Then 
$$
(f_3,f_4)=\begin{cases}
  (0,2)\\
  (1,2)\\
  (2,2)
\end{cases}
$$
The last two are forbidden by (F3), hence $f_3=0$. We are left with
$$
(f_1,f_2,f_3,f_4)=\begin{cases}
  (1,0,0,2)\\
  (1,1,0,2)\\
  (1,2,0,2)\\
  (2,0,0,2)\\
  (2,1,0,2)\\
  (2,2,0,2)
\end{cases}
$$
The configurations $(1,2,0,2), (2,1,0,2), (2,2,0,2)$ are forbidden by (F3), $(2,0,0,2)$ is forbidden by (F7), and $(1,1,0,2)$ is forbidden by {\bf Step 2} above. Therefore the only
new frequency configuration that must be excluded is $(f_1,f_2,f_3,f_4)=(1,0,0,2)$.

\end{proof}
We can summarize the result in a table.
\begin{table}[ht]
\centering
\label{tab:subpartitions}
\renewcommand{\arraystretch}{1.25}
\begin{tabular}{c|l}
\hline
Forbidden subpartition & new frequency configuration to be excluded \\
\hline
$(1,1)$ &
$
f_1=2$
\\

$(2,1)$ &
\(
(f_1,f_2)=(1,1)
\)
\\

$(2,2)$ &
\(
(f_1,f_2)=(0,2)
\)
\\

$(3,2)$ &
\(
(f_1,f_2,f_3)=(0,1,1)
\)
\\

$(3,3)$ &
\(
(f_1,f_2,f_3)=(0,0,2)
\)
\\

$(4,4,1)$ &
\(
(f_1,f_2,f_3,f_4)=(1,0,0,2)
\)
\\

\hline
\end{tabular}
\end{table}

\section*{Acknowledgements} 
\noindent
We thank Ali Uncu for making us aware of references \cite{AEH20} and \cite{Sa24}. The authors acknowledge partial support from the Royal Physiographic Society of Lund,
 the Croatian Science Foundation under the projects IP-2022-10-9006 and IP-2025-02-4720 and the project ``Implementation of cutting-edge research and its application as part of the Scientific Center of Excellence for Quantum and Complex Systems, and Representations of Lie Algebras'', Grant No. PK.1.1.10.0004, co-financed by the European Union through the European Regional Development Fund - Competitiveness and Cohesion Programme 2021-2027.

	 %%%%%%%%%%%%%%%%%%%%%%%%%%%%
	\linespread{1.0}

	\bigskip
	\noindent
	(Stefano Capparelli) Dipartimento SBAI, Universit\`a La Sapienza, Rome, Italy
    \newline
    {\em Email address:} stefano.capparelli@uniroma1.it
	
	\medskip
	\noindent
	(Arne Meurman) Department of Mathematics, Lund University, Lund, Sweden
	\newline
	{\em Email address:} arne.meurman@math.lu.se
	
	\medskip
	\noindent
	(Mirko Primc) Faculty of Science, University of Zagreb, Zagreb, Croatia
    \newline
    {\em Email address:} primc@math.hr
            
\end{document}